\documentclass[11pt,reqno,oneside]{article}
\usepackage{amsmath}
\usepackage{amssymb}
\usepackage[lmargin=1in, rmargin=0.75in, tmargin=1in, bmargin=1in]{geometry}
\usepackage{xcolor}
\usepackage{graphicx}

\usepackage{pstricks}

\newtheorem{prop}{Proposition}

\definecolor{myblue}{rgb}{0,0.1,0.9}

\definecolor{myred}{rgb}{.95,.79,0.0}

\usepackage[normalem]{ulem}

\title{On the sub-permutations of pattern avoiding permutations}
\author{Filippo Disanto\thanks{E-mail: fdisanto@uni-koeln.de}, Thomas Wiehe\thanks{E-mail: twiehe@uni-koeln.de}}
\date{\today}

\begin{document}

\maketitle

\abstract{There is a deep connection between permutations and trees. Certain sub-structures of permutations, called \emph{sub-permutations}, bijectively map to sub-trees of binary increasing trees.
This opens a powerful tool set to study enumerative and probabilistic properties of sub-permutations and
to investigate the relationships between 'local' and 'global' features using the concept of pattern avoidance.

First, given a pattern $\mu$, we study how the avoidance of $\mu$ in a permutation $\pi$ affects the presence of other patterns in the sub-permutations of $\pi$. More precisely, considering patterns of length $3$, we solve instances of the following problem: given a class of permutations $\mathcal{K}$ and a pattern $\mu$, we ask for the number of permutations $\pi \in Av_n(\mu)$ whose sub-permutations in $\mathcal{K}$ satisfy certain additional constraints on their size.

Second, we study the probability for a generic pattern to be contained in a random permutation $\pi$ of size $n$ without being present in the sub-permutations of $\pi$ \emph{generated} by the entry $1 \leq k \leq n$. These theoretical results can be useful to define efficient randomized pattern-search procedures based on classical algorithms of pattern-recognition, while the general problem of pattern-search is NP-complete.
}

\section{Introduction}

Characterizing a class of objects by means of its sub-structures is a general theme in mathematics. Using the concept of {\it pattern avoidance}, 
we investigate sub-structures of permutations, called \emph{sub-permutations}, and demonstrate
bijective relations and enumerative and probabilistic results.

Sub-permutations, as defined in Section~\ref{intro}, correspond to the classical concept of sub-trees via a well known \cite{goulden, stan} bijection which maps the set of pemutations of size $n$ onto the set of binary \emph{increasing} trees (BITs) with $n$ nodes.

The study of patterns in sub-permutations of un-constrained random permutations has been initiated \cite{inria} by Flajolet et al. Using  the equivalent sub-tree teminology for BITs, they focus on the number of times a fixed permutation occurs as a sub-permutation of a larger random permutation.

There is a long standing  interest in studying discrete structures constrained by pattern conditions. The concept of pattern avoidance was introduced by MacMahon \cite{mac} and then elaborated by Knuth in \cite{knu} and in \cite{simion} by Simion and Schmidt. Lately, several combinatorial structures -- trees \cite{dai,inria,flajosipala,row}, lattice paths \cite{berni,sapo} and compositions of integers \cite{heu} --  have been analyzed using patterns constraints.

Here, we expand these concepts to sub-permutations and, by means of equivalence, to sub-trees of BITs, which are constrained by pattern avoiding conditions. This new point of view is introduced with the general aim of better understanding the relationships  between 'local' and 'global' features of pattern avoidance.

In Section~\ref{ottobre} we concentrate on patterns of length three. We derive several enumerative and bijective properties for sub-permutations, relating them with other well-studied combinatorial structures such as planar binary trees and lattice paths. These properties are most easily formulated in the light of the following  general problem. 

\smallskip

\noindent \textbf{Problem:} given a class of permutations $\mathcal{K}$, we ask for the number of permutations $\pi \in Av_n(\mu)$  whose sub-permutations in $\mathcal{K}$ satisfy certain additional constraints on their size.
 
\smallskip

In particular, we consider $\mu=312$ and $\mu=123$. These two patterns can be taken as representative of the two clusters of patterns $\{312,213,132,231  \}$ and $\{ 123,321 \}$ which are  not equivalent with respect to the standard operators of mirror image, complement and inverse. 
In Section~\ref{ottobre} we consider four possible instances of class $\mathcal{K}$. For $\mu=312$ we consider $\mathcal{K}=Av(213)$ and $\mathcal{K}$ equal to the set of odd alternating permutations. In particular we prove that, for $n \rightarrow \infty$, taking at random a permutation in $\pi \in Av_n(312)$ the expected value of the size of the biggest sub-permutation of $\pi$ belonging to $Av(213)$ grows with $\log_2(n)$. For $\mu=123$ we consider $\mathcal{K}=Av(21)$ and $\mathcal{K}=Av(12)$. For the latter, we show that the number of permutations of $Av_n(123)$, whose biggest non-trivial decreasing sub-permutation has size bounded by a fixed number $j$, is related to the number of \emph{Dyck} paths of size $n$ avoiding the pattern $U^{j+2}D$, where $U$ (resp. $D$) stand for up (resp. down) step in the path.
 
In Section~\ref{fw} we introduce the concept of pattern avoidance in terms of sub-permutations. We find the probability to detect a pattern $213$ in a permutation by looking just at the sub-permutation generated by the entry $k=2$. Then, we generalize our results considering an arbitrary pattern $\sigma$ and a general value of the parameter $k$ and we show that our theoretical results are in agreement with data generated by Monte Carlo random experiments. 

We expect that this new kind of pattern related problems will draw a number of practical applications. For instance, for computational tasks \cite{alberto,bosse}, it can be important to know when global  properties of a permutation can be predicted, if only local properties are known or accessible. 
The kind of questions formulated here in the context of  permutations and BITs are generic for all combinatorial objects for which a notion of \emph{sub-structure} is defined.

\section{Preliminaries}\label{intro}

The set of permutations of size $n$ is denoted by $\mathcal{S}_n$ and $\mathcal{S}=\bigcup_n \mathcal{S}_n$. We assume the reader to be familiar with the classical concept of pattern avoidance in permutations. 
The subset of $\mathcal{S}_n$ made of those permutations avoiding the pattern $\sigma$ is usually denoted by $Av_n(\sigma)$ and $Av(\sigma)=\bigcup_n (Av_n(\sigma))$.

\begin{figure}
\begin{center}
\includegraphics*[scale=.5,trim=0 0 0 0]{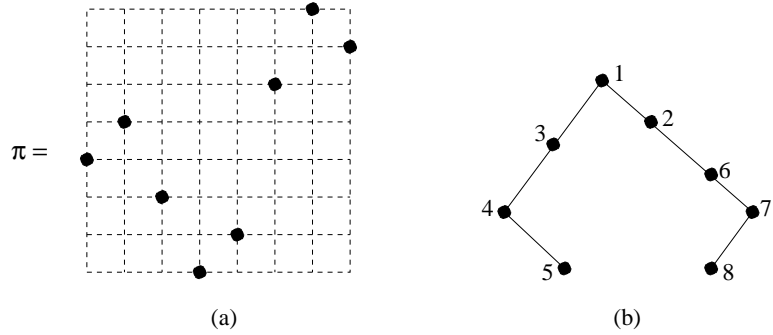}
\end{center}
\caption{(a) The permutation $\pi =  4 5 3 1 2 6 8 7$; (b) the tree associated with $\pi $ by $\phi$.}\label{permutaz}
\end{figure}

Let $\pi=\pi_1 \pi_2 \dots \pi_n$ be a permutation. For a given entry $\pi_i$ we define $s_{\pi}(\pi_i)$ as the biggest sub-string\footnote{sub-strings are sub-sequences of consecutive elements of $\pi$.} of $\pi$ which contains $\pi_i$ and whose entries are greater than or equal to $\pi_i$. Furthermore let $g_{\pi}(\pi_i)$  be the permutation 
obtained rescaling $s_{\pi}(\pi_i)$.
We call $g_{\pi}(\pi_i)$ the \emph{sub-permutation} of $\pi$ generated by $\pi_i$. 
The set of sub-permutations $(g_{\pi}(\pi_i))_{i=1\dots n}$  is  denoted by $G_{\pi}$.
Note that not all sub-strings of $\pi$ are (once rescaled) sub-permutations of $\pi$. It is crucial that only those maximal sub-strings are extracted from $\pi$ whose elements share the property of being greater than or equal to a given entry of the sub-string.
As an example consider the permutation $\pi$ which is depicted in Fig.~\ref{permutaz}~(a). In this case $G_{\pi}$ is made of 
$g_{\pi}(5)=g_{\pi}(8)=  1$, $g_{\pi}(4)  =  1 2$, $g_{\pi}(7)  =  2 1$, $g_{\pi}(3)  =  2 3 1$, $g_{\pi}(6)  =  1 3 2$,  $g_{\pi}(2)  =  1 2 4 3$ and
$g_{\pi}(1)  =  4 5 3 1 2 6 8 7$.

We also observe that there is a strong correspondance between sub-permutations as defined above and the so-called \emph{x-factorizations} of a permutation, see \cite{gabbor} and related references. Indeed, following the notation of \cite{gabbor}, the sub-permutation $g_{\pi}(\pi_i)$ is obtained by extracting together from the permutation $\pi$ the entry $\pi_i$ and the $\lambda$ and $\rho$ parts of the $\pi_i$-factorization of $\pi$.

The concept of sub-permutations is related to the one of sub-trees. To illustrate this correspondance we make use of a well known bijection \cite{goulden,stan} between the set $\mathcal{S}_n$ and the set of binary increasing trees of size $n$, denoted by $\mathcal{T}_n$. We recall that a planar rooted tree $t$, having $n$ nodes, belongs to $\mathcal{T}_n$ when:
\begin{itemize}
\item[-] each node has outdegree $0,1$ or $2$. Nodes of outdegree $0$ are called \emph{leaves};
\item[-] each node (except for the root) can be left or right oriented with respect to its direct ancestor (see Fig.~\ref{permutaz}~(b));
\item[-] each node is labelled (bijectively) with a number in $\{ 1,...,n\}$ in such a way going from the root of $t$ to any leaf of $t$ we find an increasing sequence of numbers.
\end{itemize}

The bijection $\phi:\mathcal{T}_n \rightarrow \mathcal{S}_n$ is given by the following procedure:

\begin{itemize}
\item[$i)$] given a tree $t$ each leaf of $t$ collapses into its direct ancestor whose label is then modified receiving on the left the label of the left child (if any) and on the right the label of the right child (if any). We obtain in this way a new tree whose nodes are labelled with sequences of numbers; 
\item[$ii)$] starting from the obtained tree go to step $i)$.
\end{itemize}

The algorithm $\phi$ ends when the tree $t$ is reduced to a single node whose label is then a permutation $\phi(t)$ of size $n$. For an example see Fig.~\ref{permutaz}. 

The link between sub-permutations and sub-trees is expressed by the following proposition.

\begin{prop}\label{forma}
Given a permutation $\pi = \phi(t)$, let $t_i$ be the (re-scaled) sub-tree of $t$ generated by the node $\pi_i$, then the sub-permutation $g_{\pi}(\pi_i)$ is equal to $\phi(t_i)$.
\end{prop}

By Proposition~\ref{forma}, we have that the size of the sub-permutation $g_{\pi}(\pi_i)$ is given by the number of nodes descending from node $\pi_i$ in the tree $t = \phi^{-1}(\pi)$. In \cite{kuba} the authors study statistics related to the number of descendants in subtrees of BITs. In particular, for a random permutation $\pi$ of size $n$ expectation and variance of the size of the sub-permutation generated by entry $k$ are
\begin{eqnarray}\label{attesa}
\text{E}(|g_{\pi}(k)|) &=& \frac{2n-k+1}{k+1} \,\,\, \text{and} \\
\text{Var}(|g_{\pi}(k)|) &=& \frac{2(n+1)(k-1)(n-k)}{(k+1)^2(k+2)}. \label{varianzia}
\end{eqnarray} 
More generally, for a random permutation of size $n$, the size of the sub-permutation generated by $k$ is equal to $m$ with probability
\begin{equation}\label{occhiali}
\text{Prob}(|g_{\pi}(k)| = m) = \frac{k \, m \, {{n-m-1}\choose{k-2}}}{n {{n-1}\choose{k-1}}}.
\end{equation}

\section{On the sub-permutations of $\pi \in Av(\mu)$ belonging~to~$\mathcal{K}$}\label{ottobre}

In this section we determine the number of permutations in $Av_n(\mu)$ whose sub-permutations in $\mathcal{K}$ satisfy certain size-constraints. We consider several instances of $\mu$ and $\mathcal{K}$ and find new enumerative results depending on the constraints we require each time.Some of them are connected to already known combinatorial problems.
We will study the patterns $\mu=312$ and $\mu=123$. It is well known that both the classes $Av(312)$ and $Av(123)$ are enumerated by \emph{Catalan} numbers whose first terms are $$c_0=1,c_1=1,c_2=2,c_3=5,c_4=14,c_5=42,c_6=132,c_7=429.$$ The associated ordinary generating function is \cite{ancomb}
 
\begin{equation}\label{catralan}
C(x)=\sum_{n\geq 0}c_{n}x^n=\frac{1- \sqrt{1-4x}}{2x},
\end{equation}
which has a closed form for its $n$-th coefficient given by

$$c_{n}=\frac{1}{n+1} \, {{2n}\choose{n}}.$$
Furthermore, it is well known that the asymptotic behaviour for $c_n$ is

$$c_n \sim \frac{4^{n}}{\sqrt{\pi n^3}}.$$

\subsection{$\mu=312$ and binary rooted planar trees}

We denote by $\mathcal{B}_n$ the set of binary rooted planar trees with $n$ internal nodes, where each internal node has outdegree two. It is well known that one can bijectively map the set $\mathcal{B}_{n}$ onto the set $Av_{n}(312)$. In particular we use a  bijection $\psi: \mathcal{B}_{n} \rightarrow Av_{n}(312)$ which works similarly to the mapping $\phi$ of Section~\ref{intro}, as follows. 

\begin{figure}
\begin{center}
\includegraphics*[scale=.55,trim=0 0 0 0]{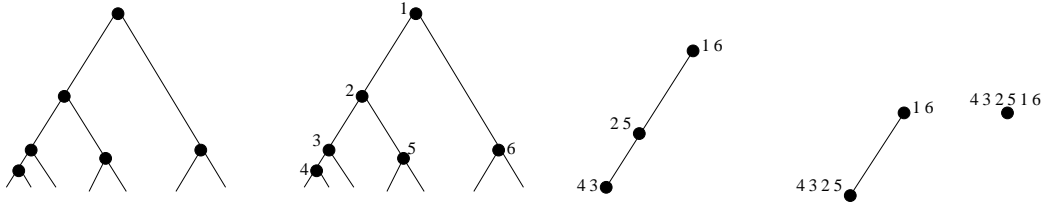}
\end{center}
\caption{The mapping $\psi$.}\label{bij}
\end{figure}

Take $t \in \mathcal{B}_{n}$ and visit its nodes according to the pre-order traversal labelling each node of outdegree two in increasing order starting with the label $1$ for the root. After this first step one has a tree labelled with integers at its nodes of outdegree two. Each leaf now collapses to its direct ancestor which takes a new label receiving on the left (resp. right) the label of its left (resp. right) child. We go on collapsing leaves until we achieve a tree made of one node which is labelled with a permutation of size $n$. See Fig.~\ref{bij} for an instance of this mapping.  

\bigskip

\subsubsection{$\mathcal{K}=Av(213)$ and caterpillar sub-trees}\label{mu1}
Let $\mathcal{K}=Av(213)$. A tree $t$ in $\mathcal{B}_n$ is a \emph{caterpillar} of size $n$, if each node of $t$ is a leaf or has  at least one leaf as a direct descendant. Through the bijection  $\psi$ we see that, given a permutation $\pi \in Av(312)$, the parameter $\gamma_{\pi}(\mathcal{K})$, that is the size of the biggest sub-permutation of $\pi$ belonging to $Av(213)$, corresponds to the size (number of internal nodes) of the biggest sub-tree of $\psi^{-1}(\pi)$ which is a caterpillar.

Let $$P_j = P_j(x)=\sum_{n \geq 0} v_{j,n} x^n$$ be the  ordinary generating function counting those permutations  $\pi$ in $Av_n(\mu)$ having $\gamma_{\pi}(\mathcal{K}) \leq j$,
we observe that $P_j$ must then satisfy the equation 

\begin{equation} \label{rec}
P_j = 1 + x{P_j}^2 - 2^{j} x^{j+1}.
\end{equation}

Indeed, a tree $t$ whose biggest caterpillar sub-tree is of size at most $j$, is either a leaf or it is built by appending to the root of $t$  two trees $t_1$ and $t_2$ of the same class. We must  exclude the case in which one of the two, $t_1$ or $t_2$, has size $0$, i.e., is a leaf, and the other one is a caterpillar of size $j$. Since there are exactly $2^{j-1}$ caterpillars of size $j$ the previous formula follows.  

From (\ref{rec}) we obtain 

\begin{equation}\label{ordinato}
P_j(x)=\frac{1-\sqrt{1 - 4x + 2^{j+2}x^{j+2}}}{2x}
\end{equation}

which gives an expansion, for $j=1$, with the following coefficients

$$1,1,0,1,2,6,16,45,126,358,1024$$ which correspond to sequence $A025266$ of \cite{sloane}. As stated there, they also count classes of Motzkin paths satisfying some constraints.

The asymptotic behaviour of the coefficients of $P_j$ follows by standard analytic methods \cite{ancomb}. When $n$ becomes large the coefficients grow like  

\begin{equation}\label{lavecchia}
v_{j,n}\sim \frac{1}{4} \sqrt{\frac{4 \rho_j -(j+2)2^{j+2} \rho_j^{j+2}}{\pi n^3}} \left( \frac{1}{\rho_j} \right)^{n+1},
\end{equation}
where $\rho_j$ is the smallest positive solution of $$1 - 4x + 2^{j+2}x^{j+2}=0.$$

Starting from the following inequality $$\frac{1}{4} < \rho_j < \frac{2}{5},$$
which is strightforward to prove for all $j \geq 1$, and applying the so-called \emph{bootstrapping} technique \cite{ancomb} one can show that $\rho_j$ tends to $1/4$ as $$\rho_j = \frac{1}{4} + \frac{1}{2^{j+4}} + O\left( j  \left( \frac{2}{5}\right)^{j}  \right).$$ 

The main result we want to present in this section concerns the average of $\gamma = \gamma_{\pi}(\mathcal{K})$ when the size of the permutation $\pi$ is large.

\begin{prop}\label{pallu}
Let $\pi \in Av_n(\mu)$, the expected value of $\gamma_{\pi}(\mathcal{K})$ is asymptotically equivalent to $\log_2(n).$
\end{prop} 
\emph{Proof.}
If $n\geq 1$ we can express the desired average value as

\begin{align}\nonumber
E_n(\gamma)&=\frac{1 v_{1,n} + \sum_{j \geq 1} (j+1)(v_{j+1,n} - v_{j,n})}{c_n} \\\nonumber
&= \frac{-v_{1,n}-...-v_{n-1,n}+nv_{n,n}+ \sum_{j \geq n} (j+1)(v_{j+1,n} - v_{j,n})}{c_n} \\\nonumber
&= \frac{-v_{1,n}-...-v_{n-1,n}+nc_n+ \sum_{j \geq n} (c_n-v_{j,n})}{c_n} \\\nonumber
&= \frac{\sum^{n-1}_{j=1} (c_n-v_{j,n}) +c_n+\sum_{j \geq n} (c_n-v_{j,n})}{c_n} \\\nonumber
&= \frac{c_n+ \sum_{j \geq 1} (c_n-v_{j,n})}{c_n} \\\nonumber
&= 1 + \frac{\sum_{j \geq 1} (c_n-v_{j,n})}{c_n} \nonumber
\end{align}

In the previous calculation we have used the fact that for $j \geq n$ we always have $v_{j,n}=c_n$.

It is sufficient now to find the $n$-th term of the function $$U(x)=\sum_{j \geq 1} (C(x)-P_j(x))= \frac{\sqrt{1-4x}}{2x} \sum_{j \geq 1} \left( \sqrt{1+\frac{2^{j+2} x^{j+2}}{1-4x}} -1 \right).$$

In what follows we want to find a function $\tilde{U}$ which estimates $U$ near the dominant singularity $1/4$. According to \cite{ancomb}, when $n$ is large, the $n$-th term of the Taylor expansion of $\tilde{U}$ provides an approximation of $[x^n]U(x)$. Our approach results to be similar to the one used in Section 3 of \cite{flajosipala}.

Let us fix $x$ near $1/4$ and let us consider the threshold function $$j_0 = \log_2 \frac{1}{|1-4x|}.$$ Then, supposing $j \geq j_0-1$, we have that $$\sqrt{1+ \frac{2^{j+2} x^{j+2}}{1-4x}} \sim \sqrt{1+ \frac{1}{2^{j+2} (1-4x)}} \sim 1 + \frac{1}{2^{j+3} (1-4x)},$$

while if  $j< j_0-1$ we can use the approximation 
$$\sqrt{1+ \frac{2^{j+2} x^{j+2}}{1-4x}} \sim \sqrt{1+ \frac{1}{2^{j+2} (1-4x)}} \sim \sqrt{ \frac{1}{2^{j+2} (1-4x)}}.$$

For $x$ sufficiently close to $1/4$  we estimate $U(x)$ as 

\begin{align} \nonumber
U(x)&\sim  \frac{\sqrt{1-4x}}{2x \sqrt{1-4x}}\left( \sum_{j \geq 1}^{j_0-2} \sqrt{\frac{1}{2^{j+2}}}\right) - \frac{\sqrt{1-4x}}{2x}\left( \sum_{j \geq 1}^{j_0-2} 1\right) + \frac{\sqrt{1-4x}}{2x(1-4x)} \left(\sum_{j \geq j_0-1} \frac{1}{2^{j+3}}\right) \\\nonumber
&=\frac{1}{4x } \sum_{j \geq 1}^{j_0-2} \sqrt{\frac{1}{2^{j}}} -\frac{\sqrt{1-4x}}{2x}(j_0-2) + \frac{1}{16x\sqrt{1-4x}} \sum_{j \geq j_0-1} \frac{1}{2^{j}} \\\nonumber
&= \frac{1}{4x} \frac{-\sqrt{2} + 2\sqrt{2}\sqrt{2^{-j_0}}}{-2 +\sqrt{2}} -\frac{\sqrt{1-4x}}{2x} \left( \log_2 \left( \frac{1}{|1-4x|} \right)-2 \right) 
 + \frac{2^{2-j_0}}{16x\sqrt{1-4x}} \\\nonumber
&=\frac{\log(2)}{x \log(16)}+\frac{\sqrt{2} \log(2)}{x \log(16)}-\frac{2 \sqrt{2-8 x} \log(2)}{x \log(16)} -\frac{\sqrt{1-4 x} \log(2)}{x \log(16)}-\frac{2 \sqrt{1-4 x} \log \left(\frac{1}{4-16 x}\right)}{x \log(16)} \\\nonumber
& \sim -\frac{2 \sqrt{1-4 x} \log \left(\frac{1}{1-4 x}\right)}{x \log(16)}. \\\nonumber
\end{align}

Using the previous calculation we can say that

\begin{equation} \label{palla}
\tilde{U}(x)=-\frac{2 \sqrt{1-4 x} \log \left(\frac{1}{1-4 x}\right)}{x \log(16)}
\end{equation}
approximates $U(x)$ near its dominant singularity $1/4$. It follows that when $n  \rightarrow \infty$

$$E_n(\gamma) \sim  \frac{[x^n]\tilde{U}(x)}{c_n}. $$

Applying standard methods \cite{ancomb,flajosipala} to (\ref{palla}) we find that

$$[x^n]\tilde{U}(x) \sim \frac{4^{n+1} \log(n)}{\log(16)\sqrt{\pi n^3}}.$$ Dividing by the asymptotic behaviour of Catalan numbers gives the claim. $\Box$

As a test one can consider the following table where, for several values of $n$, we compare the true $E_n(\gamma)$ with the approximation given by Proposition~\ref{pallu}.

\begin{center} 
\begin{tabular}{|c|ccccccc|}
\hline
$n$ & 10 & 20 & 50 & 100 & 200 & 500 & 1000 \\\hline
$E_n(\gamma)$ & 3.596 & 4.172 & 5.227 & 6.121 & 7.058 & 8.336 & 9.319 \\ 
$\log_2(n)$ & 3.321 & 4.321 & 5.643 & 6.643 & 7.643 & 8.965 & 9.965 \\\hline
\end{tabular} 
\end{center}

\bigskip

To conclude this section it is interesting to observe another possible application of the function $P_j$ described in (\ref{ordinato}). Indeed the number of permutations in $Av_n(\mu)$ having no sub-permutation of size $j$ in $\mathcal{K}$ is given by the $n$-th coefficient of

\begin{equation} \label{racchietta}
P_{j-1}(x)=\frac{1-\sqrt{1-4x+2^{j+1}x^{j+1}}}{2x}.
\end{equation}

We will compare this result with the analogous one provided in the next section.

\bigskip

\subsubsection{When $\mathcal{K}$ is the set of (odd) alternating permutations}\label{mu2}
Let $\mathcal{K}$ be the set of (odd) alternating permutations. A tree in $\mathcal{B}_{n}$, where $n$ is odd, is called \emph{strictly binary} if, removing the leaves, the remaining nodes have either out-degree $0$ or out-degree $2$. The corresponding (through the mapping $\psi$) sub-set of $Av_n(\mu)$ consists of permutations $\pi= \pi_1 \pi_2 ... \pi_n$ characterized by the following property: either $n=0,1$ or $\pi_1>\pi_2<\pi_3>...<\pi_n \, $. It is well known that the number of strictly binary trees of size $2m + 1$ is $c_{m}$. 

The ordinary generating function $L_j = L_j(x)$, with $j=2m+1$, counts  those trees in $\mathcal{B}_n$ having \emph{at least} one strictly binary sub-tree of size $j$. Equivalently, it can be seen as the function counting the permutations in $Av_n(\mu)$ with at least one alternating sub-permutation of size $j$. $L_j$ must then satisfy

\begin{equation} \label{reccc}
L_j = c_{m}x^j + x{L_j}^2 + 2xL_j (C-L_j),
\end{equation}
where $C = C(x)$ is the generating function of Catalan numbers as in (\ref{catralan}).

Indeed a tree counted by $L_j$ is either a strictly binary tree of size $j$ (first summand) or it is built by appending to the root two trees, where at least one of them contains a strictly binary sub-tree of size $j$ (second and third summand).  

By solving (\ref{reccc}) we obtain 

\begin{equation}\label{ordinatooo}
L_j(x)=\frac{\sqrt{1-4x+4c_{m}x^{j+1}}-\sqrt{1-4x}}{2x}.
\end{equation}

From (\ref{ordinatooo}, we can also  determine the number of those permutations avoiding the pattern $312$ and without (odd) alternating sub-permutation of size $j$. This is given by 
\begin{equation}\label{oordinatooo}
C(x)-L_j(x)=\frac{1-\sqrt{1-4x+4c_{m}x^{j+1}}}{2x}.
\end{equation}

It is now possible to compare results (\ref{racchietta}) and (\ref{oordinatooo}) for a fixed $j=2m + 1$ when $m$ is large. In other words, we want to measure how the avoidance of the pattern $\mu = 312$ affects the presence of the pattern $213$ in sub-permutations (from
(\ref{racchietta})) in comparison with the avoidance of alternating sub-permutations (from (\ref{oordinatooo})). By means of the bijection $\psi$ with trees described before, this gives a comaprison between caterpillar-shaped subtrees and strictly binary subtrees. 
Let us consider $a_m > 1/4$ (resp. $b_m > 1/4$) defined as the smallest positive root of $1-4x+4c_{m}x^{2m+2}=0$ (resp. $1-4x+2^{2m+2}x^{2m+2}=0$). By asymptotic considerations, we know that when $m\rightarrow \infty$ both $a_m$ and $b_m$ tend to $1/4$. But we can be more precise: when $m$ is large, the equality
$$1-4a_m+4c_{m}a_m^{2m+2}=0=1-4b_m+2^{2m+2}b_m^{2m+2}$$ implies
$$\frac{a_m}{b_m} = \frac{-1+4^m b_m^{2m+1}}{-1+c_m a_m^{2m+1}} \sim \frac{-1+(\frac{1}{4})^{m+1}}{-1+c_m (\frac{1}{4})^{2m+1}} < 1.$$ For example, if $m=5$, we have that $a_m/b_m=0.999765$ while the previous approximation gives $0.999766$.

By standard methods \cite{ancomb} one can compute the asymptotic behaviour of the $n$-th coefficient of (\ref{oordinatooo}) as $$\frac{1}{4} \sqrt{\frac{4 a_m -(2m+2) 4c_m a_m^{2m+2}}{\pi n^3}} \left( \frac{1}{a_m} \right)^{n+1}  $$ while the behaviour of the $n$-th coefficient of (\ref{racchietta}) is given by (\ref{lavecchia}) (considering $b_m=\rho_{j-1}$) as $$\frac{1}{4} \sqrt{\frac{4 b_m -(2m+2)2^{2m+2} b_m^{2m+2}}{\pi n^3}} \left( \frac{1}{b_m} \right)^{n+1}.$$ From these results follows 

\begin{prop}\label{rain}
For a fixed and sufficiently large $m$, when $n\rightarrow \infty$, the ratio
\begin{equation}\label{rascio}
\frac{|\{ \pi \in Av_n(\mu): (G_{\pi}\cap Av(213)\cap \mathcal{S}_{2m+1}) = \emptyset  \}|}{|\{ \pi \in Av_n(\mu): (G_{\pi}\cap \mathcal{K}\cap \mathcal{S}_{2m+1}) = \emptyset\}|}
\end{equation}
goes to $0$ equally fast as
\begin{equation}\label{testa}
k_m \left( \frac{-1+(\frac{1}{4})^{m+1}}{-1+c_m (\frac{1}{4})^{2m+1}} \right)^{n+1},
\end{equation}
where $k_m$ is a constant depending only on $m$.
\end{prop}

In the following table we compare the values of the ratio of Proposition~\ref{rain} with the asymptotic ratio (\ref{testa}) for $m=5$ and different values of $n$.

\begin{center} 
\begin{tabular}{|c|ccccc|}
\hline
$n$ & 50 & 500 & 1000 & 5000  & 10000   \\\hline
ratio (\ref{rascio})  & 0.986  & 0.887  & 0.789  & 0.308  & 0.095    \\
$k_m=1$ in (\ref{testa}) & 0.988  & 0.889  & 0.791  & 0.310  & 0.096   \\\hline
\end{tabular} 
\end{center}

\subsection{$\mu=123$, a generating tree approach}\label{mu3}

\begin{figure}
\begin{center}
\includegraphics*[scale=0.9,trim=0 0 0 0]{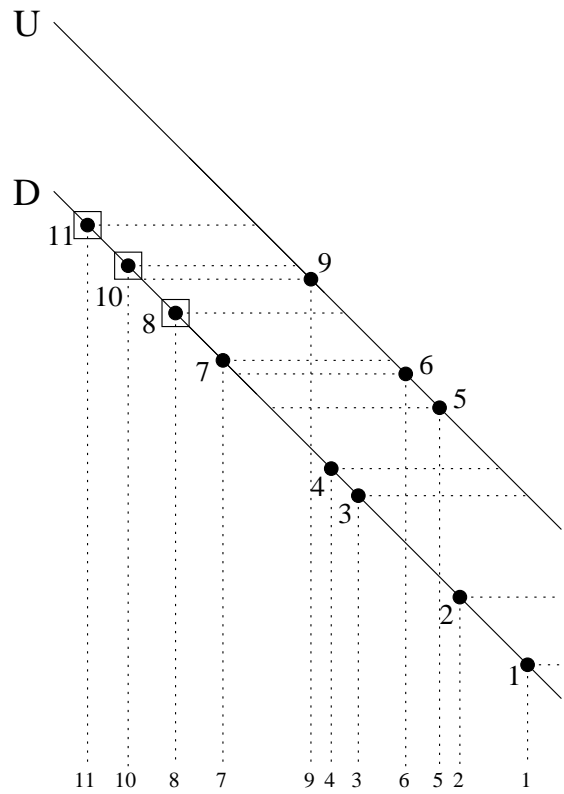}
\end{center}
\caption{The two-line representation of $\pi = 11 \, 10 \, 8 \, 7 \, 9 \, 4 \, 3 \, 6 \, 5 \, 2 \, 1 \in Av_{11}(123)$; for example $6$ covers $4$  and $3$. The entries that generate trivial decreasing sub-permutations are in boxes.}\label{deco2}
\end{figure}

In what follows let $\mu=123$ be fixed. Permutations which are the union of two decreasing subsequences have been studied since long \cite{knu}. It is known that they are characterized by the property of having no increasing subsequence of length $3$. Thus the entries of a permutation $\pi \in Av(\mu)$ can be seen as points lying on two non-intersecting lines as depicted in Fig.~\ref{deco2}. Viceversa, each permutation which can be drawn in such a way avoids $\mu$. A point $p_1$ on the line $D$ is  \emph{covered} by a point $p_2$ of the upper line $U$ if $p_2$ is on the right and above $p_1$. In order to avoid redundancies in the two-line representation of a permutation we have to respect the following rule: a point belongs to the upper line  $U$ if and only if it covers at least one point of $D$.

\begin{figure}
\begin{center}
\includegraphics*[scale=0.45,trim=0 0 0 0]{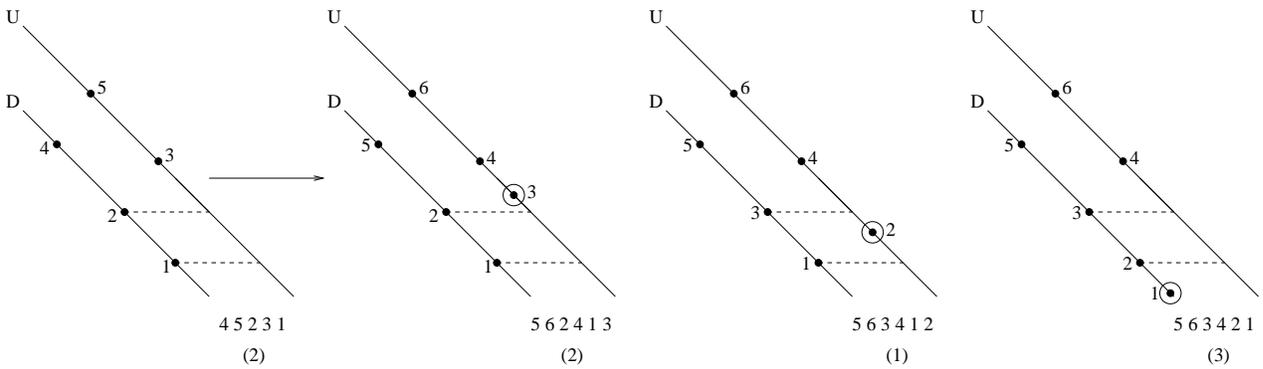}
\end{center}
\caption{The permutations produced starting from $\pi = 4 5 2 3 1$. The new rightmost entry is highlighted.}\label{passo}
\end{figure}

It is useful to observe that the set $Av_{n+1}(\mu)$ can be generated by the permutations of $Av_n(\mu)$ adding the rightmost entry. Taken $\pi \in Av_n(\mu)$, let us define $u(\pi)$ as the right-most point placed on the line $U$ (if any) and let $l(\pi)$ be the number of elements which are placed on $D$ with a smaller ordinate than $u(\pi)$. In order to create a permutation $\pi'$ of size $n+1$ we add on the right of $\pi$ a new element $p'$. If $p'$ is placed on the $D$ line then $l(\pi')=l(\pi)+1$. Otherwise $p'$ can be placed on the $U$ line in exactly $l(\pi)$ different ways, see Fig.~\ref{passo}. In this case $l(\pi')$ ranges between $1$ and $l(\pi)$. 

In conclusion, we have described a \emph{generating tree} procedure which starts with the single element permutation and, always adding the rightmost entry, generates all the permutations in $Av(123)$. Furthermore, if we consider $(l)=(l(\pi))$ the previous procedure can be summarized as

\begin{equation}\label{regole}
(l) \rightsquigarrow (1),(2),...,(l+1). 
\end{equation}

This rule can be viewed as a tree generating procedure, where a label $(l)$ is a node with $(l+1)$ descendants. If the starting label (i.e. the root of the tree) is $(0)$, the number of nodes at level $n$ is given by $c_n$. More details on generating trees and their applications to enumerative problems can be found in \cite{BB_MDF,eco,west}.   

\bigskip

\subsubsection{$\mathcal{K}=Av(21)$, the biggest increasing sub-permutation}
Let $\mathcal{K} = Av(21)$ and let $\gamma_{\pi}(\mathcal{K}) \in \{  1,2 \}$ be the size of biggest increasing sub-permutation of a given $\pi \in Av(\mu)$. Notice that $\gamma_{\pi}(\mathcal{K})=2$ if and only if the first (from left to right) three entries of $\pi$ are in the order $231$. Indeed, if $a<b$ are the entries of an increasing sub-permutation of size $2$, such  a sub-permutation is the one generated by $a$. It follows that $a$ must be the leftmost entry of $\pi$ while $b$ the topmost. In order to 'close' the sub-permutation, we also need the presence of a point $c<a$ placed immediately on the right of $b$.
For instance, the permutation $\pi = 132$ has its biggest increasing sub-permutation of size $1$ (and not $2$) because, with notations of Section~\ref{intro}, one has $g_{\pi}(1) = 132, g_{\pi}(2)=21$ and $g_{\pi}(3) = 1$.

To count the permutations $\pi \in Av(\mu)$ of size $n(\pi)\geq 3$ with the biggest increasing sub-permutation of size $2$ we have just to consider rule (\ref{regole}) with a starting point given by $(2)=l(231)$. 

The ordinary generating function $$M_2(x,y)=\sum_{\pi}x^{n(\pi)}y^{l(\pi)}$$ satisfies the functional equation
\begin{equation}\label{orabasta}
\left( 1+\frac{xy^2}{1-y} \right) M_2(x,y) = x^3y^2 + \frac{xy}{1-y} M_2(x,1) 
\end{equation}

which can be solved by the \emph{kernel method} \cite{BB_MDF,prodi} finding

$$M_2(x,1)=\frac{1}{4} \left(-1+\sqrt{1-4 x}\right) \left(-1+\sqrt{1-4 x}+2 x\right).$$

A closed formula for the coefficients is given by

$$[x^n]M_2(x,1)=\frac{3 (2n-4)!}{(n-3)! n!} \mathrm{\, \, (with \,} n\geq 3). $$

In conclusion we can state that

\begin{prop}
If $n\geq 3$, the number of permutations in $Av_n(123)$ having an increasing sub-permutation of size $2$ is $$a_n = \frac{3 (2n-4))!}{(n-3)! n!}.$$ Furthermore, the ones having the biggest increasing sub-permutation of size $1$ are (with $n \geq 3$) $$b_n = c_n - a_n,$$ where $c_n$ is the $n$-th Catalan number.  
\end{prop}

For $n=3,...,10$ the following table shows the values of $a_n$ and $b_n$.

\begin{center} 
\begin{tabular}{|c|cccccccc|}
\hline
$n =$ & 3 & 4 & 5 & 6 & 7 & 8 & 9 & 10   \\\hline
$a_n$ & 1 & 3 & 9 & 28  & 90 & 297  & 1001  & 3432   \\
$b_n$ & 4 & 11 & 33 & 104 & 339 & 1133 & 3861  & 13364   \\ \hline
\end{tabular} 
\end{center}

The coefficients $a_n$ are the (shifted) entries of sequence $A000245$ of \cite{sloane}, while the numbers $b_n$ do not appear there. 

\bigskip

\subsubsection{$\mathcal{K}=Av(12)$ and pattern avoidance in Dyck paths}\label{mu4}
Let $\mathcal{K} = Av(12)$. Thus, in this section, we investigate decreasing sub-permutations of a given $\pi \in Av(\mu)$. 

Based on the two line representation of $\pi$, which is described at the beginning of Section~\ref{mu3}, we can observe that there are two kinds of decreasing sub-pemutations. More precisely, if $g_{\pi}(\pi_i)$ is a decreasing sub-permutation of size $j$, either $g_{\pi}(\pi_i)$ corresponds to a sequence of $j$ adjacent points (adjacent with respect to their abscissas) lying on the line $U$ or $g_{\pi}(\pi_i)$ consists of the first $j$ entries of $\pi$ (that are therefore placed on the line $D$). In the latter case, the sub-permutations is said to be \emph{trivial} (see Fig.~\ref{deco2}).
Indeed, if $\pi_i$ is on the line $U$, then all the entries of $g_{\pi}(\pi_i)$ belong to the same line. On the other hand, if $\pi_i$ is placed on $D$, then
%, first, there is no entry of $g_i$ that is at the same time placed on $U$ and on the right of $\pi_i$ and, second, 
there is no entry of $U$ on the left of $\pi_i$. Such an entry should indeed cover an entry placed on $D$. The latter would be also present in $g_{\pi}(\pi_i)$ that would then be not decreasing.
%see that $g_{\pi}(\pi_i)$ is decreasing if and only if it corresponds to a sequence of adjacent points (adjacent with respect to their abscissas) lying on the line $U$, with the only exception of $\pi \in Av(12)$. 

It can also be noticed that the entries of a non-trivial decreasing sub-permutations are always introduced (reading from left to right) by an ascending step of the permutation.
  
For instance, the permutation depicted in Fig.~\ref{deco2} has three non-trivial decreasing sub-permutations, the one generated by the entry $9$ and introduced by the ascent $7-9$, the one generated by $6$ introduced by $3-6$ and the biggest one which is generated by $5$ and still introduced by the ascent step $3-6$.

In what follows, we denote by $\gamma^U_{\pi}$ the size of the biggest non-trivial decreasing sub-permutation of $\pi$ and we study the number of $\pi$'s such that $\gamma^U_{\pi}$ is \emph{at most} $j$.

Given $\pi$, it is useful to consider the point $d(\pi)$ defined as the right-most point placed on the line $D$. With the terminology of Section~\ref{mu3}, the parameter $v=v(\pi)$ is defined as the number of points placed on the line $U$ which are covering $d(\pi)$. In the permutation of Fig.~\ref{deco2} we have $v(\pi)=0$.
Looking at the recursive contruction of $Av(\mu)$, which has been defined in Section~\ref{mu3}, we observe that, in order to create only permutations satisfying $\gamma^U \leq j$, it is sufficient to avoid at each step the construction of a permutation with a $v$-value greater than $j$. This corresponds to use the rules of the form 
\begin{equation}\label{regolaproibita}
(l)\rightsquigarrow (l') \,\, \mathrm{(with}\,\,l' \leq l) 
\end{equation}
no more than $j$ times consecutively. 

If we put the same restrictions on a particular recursive construction holding for \emph{Dyck} paths \cite{deu}, it turns out that the statistic $\gamma^U_{\pi} \leq j$ is equivalent to determine the number of those paths avoiding the pattern $U^{j+2}D$ (defined below) and having a fixed size. This problem on paths has been deeply studied in the literature, see for example \cite{sapo} as well as the following sequences of \cite{sloane}: $A001006$, i.e., the one of \emph{Motzkin} numbers, for the case $j=1$, $A036765$ for the case $j=2$ and $A036766$ for $j=3$.

We assume the reader is familiar with the definition of Dyck paths. The set of Dyck paths of semi-length $n$ is denoted by $\mathcal{D}_n$ and it is a well known fact that $|\mathcal{D}_n|=c_n$. We represent  a path $p$  as a sequence of up $U$ and down  $D$ steps. We define a \emph{block} (or a \emph{primitive} Dyck path \cite{deu}) of $p$ as a minimal sub-string (made of consecutive entries) which is still Dyck path and which starts at height $0$ in $p$. For example, if $p=UUDUDDUD$, then there are only two blocks inside $p$: $UUDUDD$  made by the first six steps and $UD$ which corresponds to the last two entries of $p$.

The set $\mathcal{D}_{n+1}$ can be generated recursively by the paths in $\mathcal{D}_n$ according to a construction which we use to show the main result of this section. Let 
$$p = b_1 ... b_{k-1} b_k $$ 
be a path of  $\mathcal{D}_n$ decomposed in terms of its blocks, where $k>0$.
In order to create a path $p'$ of semi-length $n+1$ we add two steps $U$ and $D$ inside $p$ in $k+1$ possible ways. Obtaining 

\begin{align}\label{merlino}
p'=&  UD p \,\, \mathrm{or} \\\label{merlimo} 
p'=&  U b_1 ... b_{i} D b_{i+1} ... b_k  \,\, \,\, (1\leq i \leq k).  
\end{align}

In both cases we add an $U$ step at the beginning creating the left-most block of $p'$. In \cite{deu} this decomposition is called \emph{first return} decomposition.
If $p$ is a path, let $l(p)$ be the number of its blocks. The construction above corresponds then to the generating tree associated with the rules already described in (\ref{regole}). In particular, (\ref{merlino}) is associated with  $(l) \rightsquigarrow (l+1)$. 

Now, we are ready to prove the following result.

\begin{prop}
The number of permutations $\pi \in Av_n(\mu)$ having $\gamma^U_{\pi} \leq j$ equals the number of paths in $\mathcal{D}_n$ avoiding $U^{j+2}D$
\end{prop}
\emph{Proof.}
Considering what we have already shown previously, we have to prove two facts: $i)$ if we  use the rule $(l)\rightsquigarrow (l')$ with $l' \leq l$ for more than $j$ consecutive times in the construction of a path $p$, then $p$ contains the considered pattern; $ii)$ if $p$ contains $U^{j+2}D$, then it has been created using rule (\ref{regolaproibita}) consecutively at least $j+1$ times. To prove $i)$ it is enough to observe that we obtain a path which starts as $U^{j+1}b_1...$ if we start from a generic path $b_1 ... b_k$ and we apply (\ref{regolaproibita}) $j+1$ times consecutively. The claim follows because $b_1$ starts as $b_1= U^iD ... $, with $i>0$.
To prove $ii)$ let us suppose that $U^{j+2}D$ is inside a path $p$ and let us denote by $T$ the generating tree associated with rules (\ref{merlino}) and (\ref{merlimo}). Observe that we can move inside $T$ going from $p$ to the root simply removing, step-by-step, the left-most entry $U$ of $p$ and the corresponding $D$. We must find, at some point, an ancestor of $p$ which starts as $U^{j+2}D ...$. To buid this ancestor, the last $j+1$ applications of the construction rule belong to case (\ref{merlimo}). $\Box$        

\section{Pattern avoidance looking at sub-permutations}\label{fw}

Using the terminology of Section~\ref{intro} it is obvious that, for any pattern~$\sigma$, each permutation $\pi$ satisfies the following equivalence $$\pi \in Av(\sigma) \Longleftrightarrow g_{\pi}(1) \in Av(\sigma).$$ This is due to the fact that $g_{\pi}(1)=\pi$. This does in general not hold for an entry $k>1$. Given $k>1$ and a pattern $\sigma$, it seems interesting to consider permutations $\pi$ for which $$\pi \in Av(\sigma) \Longleftrightarrow g_{\pi}(k) \in Av(\sigma).$$ In other words we seek  $\pi$ such that if it contains $\sigma$ ($\sigma \prec \pi$), then also $\sigma \prec g_{\pi}(k)$. We denote the set of such permutations by $Av(\sigma;k)$. Furthermore note that in terms of probability one has 
$$\text{Prob}(\pi \in Av(\sigma) | g_{\pi}(k) \in Av(\sigma))=\text{Prob}(\pi \in Av(\sigma;k))$$ or equivalently
\begin{equation}\label{joe}
\text{Prob}(\sigma \prec \pi  \text{ and } \sigma \nprec g_{\pi}(k) )=\text{Prob}(\pi \notin Av(\sigma;k))
\end{equation}
 
which quantifies the presence of the pattern $\sigma$ in $\pi$ depending on its presence in the sub-permutation generated by the entry $k$.

\begin{figure}
\begin{center}
\includegraphics*[scale=.9,trim=0 0 0 0]{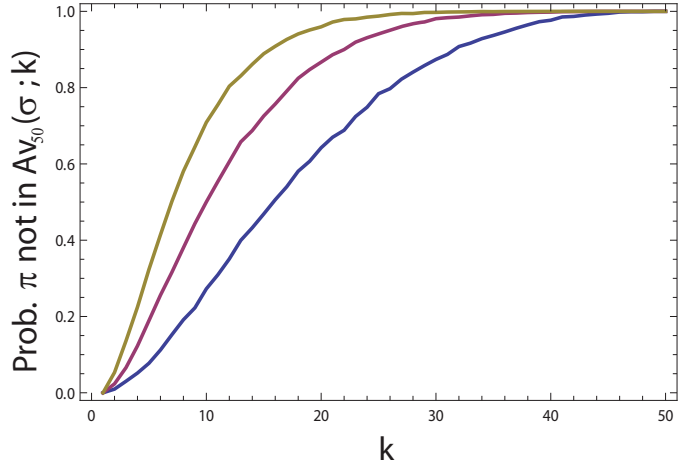}
\end{center}
\caption{Plot of the probability $\text{Prob}(\pi \notin Av_{50}(\sigma;k))$ for patterns $\sigma = 213, 1324, 25314$ (from bottom to top). }\label{vsk}
\end{figure}

For a fixed $n$, we expect the value for $\text{Prob}(\pi \notin Av_{n}(\sigma;k))$ to increase from $0$ to $1$  when $k$ goes from  $k=1$ to $k=n$. Using randomly  generated permutations, this behaviour is depicted in Fig.~\ref{vsk} for several patterns of different lengths. The size of the random pemutation $\pi$ is fixed at $n=50$.

\paragraph{Computational aspects.}
From an algorithmic point of view, it is a difficult task to decide whether a permutation $\sigma$  is contained as pattern in a text-permutation $\pi$. In fact, in its generality, this problem is NP-complete \cite{alberto,bosse}. Therefore, predicting the presence of a pattern based on the inspection of a limited part of the input text-permutation, can make an appreciable difference in terms of running-time of a pattern recognition procedure.

A sub-permutations approach and knowledge of the probability (\ref{joe}) can be used to introduce a probabilistic aspect to pattern search. Considering sub-permutations, for a 'probabilistic' version of any search algorithm one can make the following statement: in case the pattern $\sigma$ is not detected in $g_{\pi}(k)$ according to the given search procedure, then $\sigma$ is in $\pi$ with {\it probability} 
%\begin{eqnarray}\nonumber
%\text{Prob}(\sigma \in \pi | \sigma \notin g_{\pi}(k)) &=&  \frac{\text{Prob}(\pi \notin Av(\sigma;k))}{ \text{Prob}(g_{\pi}(k) \in Av(\sigma))} = \frac{\text{Prob}(\pi \notin Av(\sigma;k))}{\text{Prob}(\pi \notin Av(\sigma)) + \text{Prob}(\pi \notin Av(\sigma;k))} \\\nonumber
%&\simeq & \frac{\text{Prob}(\pi \notin Av(\sigma;k))}{1 + \text{Prob}(\pi \notin Av(\sigma;k))}. 
%\end{eqnarray}
\begin{equation}\label{pancione}
\text{Prob}(\sigma \prec \pi | \sigma \nprec g_{\pi}(k)) =  \frac{\text{Prob}(\pi \notin Av(\sigma;k))}{ \text{Prob}(g_{\pi}(k) \in Av(\sigma))}.
\end{equation}

Observe that focusing on $g_{\pi}(k)$ reduces the size of the input text-permutation for the search-procedure according to (\ref{attesa}) and (\ref{varianzia}). For instance, if $k=3$, on average we pass from an input of size $n$ to roughly $n/2$. Vice versa, if $k$ is large, say $n = \mathcal{O}(k)$, then the input for the original search procedure is on average a constant. For example, $k=n/2$ gives $E(g_{\pi}(k))) \simeq 3$. In this cases, when the expectation for $|g_{\pi}(k)|$ is a constant $K$ and assuming the knowledge of the cardinality $|Av_K(\sigma)|$, the denominator of (\ref{pancione}) can be approximated by considering $g_{\pi}(k)$ as a random permutation of size $K$. 

Below, we investigate the behaviour of $\text{Prob}(\pi \notin Av_{n}(\sigma;k))$ for random pemutations of size $n$. We start by providing an exact result for the Catalan pattern $\sigma = 213$ and the sub-permutations generated by the entry $k=2$. 

\subsection{The cardinality of $Av_n(213;2)$ and generalizations}\label{fw2}
We proceed counting those permutations $\pi \in \mathcal{S}_n $  such that $213 \prec \pi$ and $g_{\pi}(2) \in Av(213)$. Given $\pi$ of size $n$, we denote by $m$ the lowest entry of $\pi$ such that $m \neq 1$ and $m \notin s_{\pi}(2)$. Observe that $m$ could not exist. 
In order to have  $\pi \notin Av_n(213;2)$ only three distinct situations are possible.

\begin{itemize}
\item[$i)$] The entry $2$ is placed on the left of $1$ and $m$ exists;
\item[$ii)$] the entry $2$ is placed on the right of $1$ and $213 \prec g_{\pi}(m)$;
\item[$iii)$] the entry $2$ is placed on the right of $1$ and $g_{\pi}(m) \in Av(213)$. Observe that in this case we find $213 \prec \pi$ if and only if $m < M$, where $M$ is defined as the biggest entry of $\pi$ belonging to $g_{\pi}(2)$.  
\end{itemize}

Let us denote by $i$ the cardinality of $g_{\pi}(2)$, then we have $1\leq i \leq n-1$ and $n-i-1$ gives the number of entries in $g_{\pi}(m)$. 
The number of permutations corresponding to $i)$ is then
\begin{equation}\label{sco}
\sum_{i=1}^{n-2} c_i \, (n-i-1)! \, {{n-2}\choose{i-1}}.
\end{equation}
Similarly, the number of possible instances of $ii)$ is given by 
\begin{equation}\label{reg}
\sum_{i=1}^{n-4} c_i \, ((n-i-1)! - c_{n-i-1}) \, {{n-2}\choose{i-1}}.
\end{equation}
Finally, in case $iii)$ we have 
\begin{equation}\label{gia}
\sum_{i=1}^{n-2} c_i \, c_{n-i-1} \, \left({{n-2}\choose{i-1}} - 1\right). 
\end{equation}
These three results gives the following

\begin{prop}
The number of permutations of size $n$ wich are not in $Av_n(213;2)$ is given by 
\begin{align}\label{peperone}
|\mathcal{S}_n\setminus Av(213;2)| =& \,\, 2 (n-2)! \left( \sum_{i=1}^{n-4} \frac{c_i}{(i-1)!} \right) + 2 (n-2)(n-3) c_{n-3} \\\nonumber
& \,\, + 2 (n-2) c_{n-2} - c_n + 2 c_{n-1}, \nonumber
\end{align}
where $c_n$is the $n$-th Catalan number.
\end{prop}

It follows that the cardinality of $Av_n(213;2)$ is given by

$$|Av_n(213;2)|=n! - |\mathcal{S}_n\setminus Av(213;2)|$$ and, for $3\leq n \leq 10$, the numbers are $$5,16,68,392,2905,25508,251188,2703440,$$ which are not listed in \cite{sloane}.

When $n$ is large we can approximate the sum inside parentheses in (\ref{peperone}) with a constant
\begin{equation}\label{costantino}
h = \lim_{n \rightarrow \infty} \sum_{i=1}^{n-4} \frac{c_i}{(i-1)!} = 11.75330... \,\, .
\end{equation}
An asymptotic approximation of $|\mathcal{S}_n\setminus Av(213;2)|$ is then given by
$$2 (n-2)! h.$$
The probability that a permutation of size $n\rightarrow \infty$ does not belong to $Av(213;2)$ is then 
\begin{equation}\label{piupiu}
\text{Prob}(\pi \notin Av_{n}(213;2)) \sim \frac{2 h}{n^2}. 
\end{equation}

\paragraph{Generalazing to other patterns $\sigma$.}

Observe that equation (\ref{piupiu}) depends only on the constant $h$ (besides $n$). This constant can be determined with high accuracy from the first few, say $50$, summands in (\ref{costantino}). This suggests that often one does not need to know for every $n$ the number of permutations in $Av_n(\sigma)$ to determine the behaviour of $Av_n(\sigma;k)/n!$. In fact, if $k=2$ and for $n$ large,    

\begin{equation}\label{piupiupiu}
\text{Prob}(\pi \notin Av_n(\sigma;2)) \sim \frac{2 h_{\sigma}}{n^2},
\end{equation}
where $$h_{\sigma} = \lim_{n \rightarrow \infty} \sum_{i=1}^{n-1-|\sigma|} \frac{|Av_i(\sigma)|}{(i-1)!}.$$
Note that the sum in (\ref{piupiupiu}) converges because of the Stanley-Wilf bound \cite{marcus}. The constant $h_{\sigma}$ can be approximated using the first values of $(|Av_n(\sigma)|)_n$.
The reasoning which leads to (\ref{piupiupiu}) is the following. As done for the case $Av(213;2)$, let us consider a permutation $\pi \notin Av_n(\sigma;2)$. Let $i$ be the size of $g_{\pi}(2)$ and $m$ the lowest entry of $\pi$ with $m \neq 1$ and $m \notin s_{\pi}(2)$. There are two basic cases for $\pi$ depending on the presence of the pattern $\sigma$ in $g_{\pi}(m)$. If $g_{\pi}(m) \notin Av(\sigma)$, the number of possible $\pi$'s is given by 
$$2(n-2)! \left( \sum_{i=1}^{n-1-|\sigma|} \frac{|Av_i(\sigma)|}{(i-1)!} - \sum_{i=1}^{n-1-|\sigma|} \frac{|Av_i(\sigma)|}{(i-1)!} \times \frac{|Av_{n-i-1}(\sigma)|}{(n-i-1)!}  \right).$$

If $g_{\pi}(m) \in Av(\sigma)$ then the number of possible $\pi$'s is bounded from the top by 
$$2(n-2)! \sum_{i=1}^{n-1} \frac{|Av_i(\sigma)|}{(i-1)!} \times \frac{|Av_{n-i-1}(\sigma)|}{(n-i-1)!}.$$

It follows that 

\begin{align}\label{jojo}
& 2(n-2)! \left( h_{\sigma} - \sum_{i=1}^{n-1-|\sigma|}  \frac{|Av_i(\sigma)|}{(i-1)!} \times \frac{|Av_{n-i-1}(\sigma)|}{(n-i-1)!}  \right) & \\\nonumber
& \leq |\mathcal{S}_n \setminus Av(\sigma;2)| & \\\label{coco}
& \leq 2(n-2)! \left( h_{\sigma} + \sum_{i=n-|\sigma|}^{n-1} \frac{|Av_i(\sigma)|}{(i-1)!} \times \frac{|Av_{n-i-1}(\sigma)|}{(n-i-1)!} \right). &\\\nonumber
\end{align}

The sums appearing in (\ref{jojo}) and (\ref{coco}) go to zero for $n$ large. In particular, using the Stanley-Wilf bound \cite{marcus} with $|Av_n(\sigma)| \leq c^n$, the sum in (\ref{jojo}) satisfies
$$ \sum_{i=1}^{n-1-|\sigma|}  \frac{|Av_i(\sigma)|}{(i-1)!} \times \frac{|Av_{n-i-1}(\sigma)|}{(n-i-1)!} \leq \frac{ c^n \text{Pol.}(n) }{ n! } \sum_{i=0}^n {{n}\choose{i}} = \frac{ (2c)^n \text{Pol.}(n) }{ n! } \rightarrow 0.$$

Dividing by $n!$ we obtain the desired result 

$$\text{Prob}(\pi \notin Av_{n}(\sigma;2)) \sim \frac{2 h_{\sigma}}{n^2}.$$

\paragraph{Approximations for higher values of $k$.}

\begin{figure}
\begin{center}
\includegraphics*[scale=.85,trim=0 0 0 0]{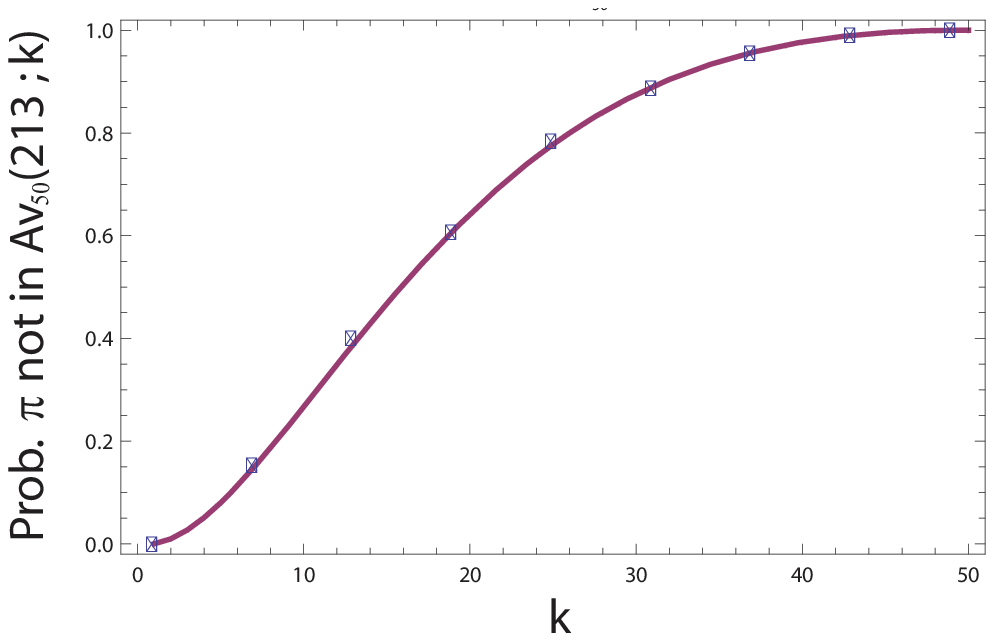}
\includegraphics*[scale=.85,trim=0 0 0 0]{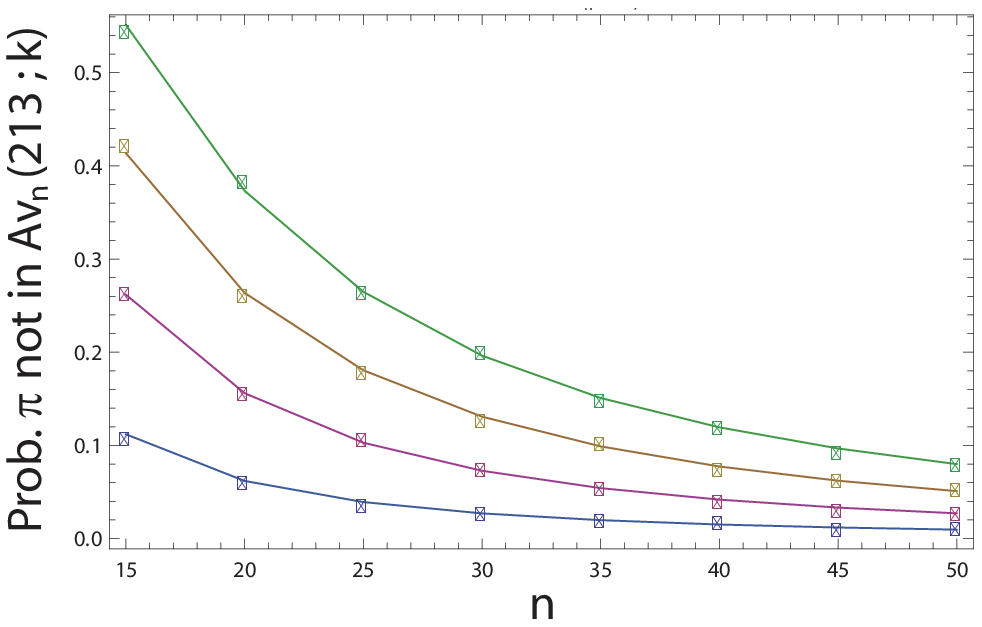}
\end{center}
\caption{(Top) Plot of $\text{Prob}(\pi \notin Av_{50}(213;k))$ for $1\leq k \leq 50$. Dots are obtained experimentally via computer simulation. The line is drawn according to (\ref{puffo}). (Bottom) Plot of $\text{Prob}(\pi \notin Av_{n}(213;k))$ for small values of $k$ and $15\leq n \leq 50$. From bottom to top we set $k=2,3,4,5$ respectively. Points correspond to experimental data while lines represent theoretical results according to (\ref{puffo}).}\label{rombodituono}
\end{figure}

If $\pi$ is a random permutation of size $n$ large enough, we can assume that
$$\text{Prob}(\sigma \nprec g_{\pi}(k) \text{ and } \sigma \prec \pi) \simeq \text{Prob}(\sigma \nprec g_{\pi}(k)).$$
Thus 
\begin{eqnarray}\nonumber
\text{Prob}(\pi \notin Av_n(\sigma;k)) &\simeq & \sum_{m=1}^{n-k+1}  \text{Prob}(|g_{\pi}(k)| = m \text{ and } \sigma \nprec g_{\pi}(k))  \\\nonumber 
&=& \sum_{m=1}^{n-k+1} \text{Prob}(|g_{\pi}(k)| = m) \times \text{Prob}(\sigma \nprec g_{\pi}(k) | \text{ size(}g_{\pi}(k)) = m).  
\end{eqnarray}
Averaging over all possible permutations $\pi$ of size $n$, we take the sub-permutation $g_{\pi}(k)$ as a random permutation of size $1\leq m \leq n-k+1$ with probability given by (\ref{occhiali}). Thus we estimate
\begin{equation}\label{puffo}
\text{Prob}(\pi \notin Av_n(\sigma;k)) \simeq  \frac{k}{n {{n-1}\choose{k-1}}}  \sum_{m=1}^{n-k+1}  \frac{ m \, {{n-m-1}\choose{k-2}}}{m!}  \cdot |Av_m(\sigma)|.
\end{equation}
For instance, if $|\sigma| = 3$, (\ref{puffo}) can be computed using Catalan numbers
and the result can be compared with simulation-data used to obtain Fig.~\ref{vsk} (bottom line). This is done in Fig.~\ref{rombodituono}.

Note that (\ref{puffo}) can be rewritten as
\begin{equation}\label{buffetto}
\text{Prob}(\pi \notin Av_n(\sigma;k)) \simeq \frac{k(k-1)}{n(n-1)} \sum_{m=1}^{n-k+1} \frac{{{n-m-1}\choose{k-2}}}{{{n-2}\choose{k-2}}} \times \frac{|Av_m(\sigma)|}{(m-1)!}.
\end{equation}
This shows that, from a theoretical point of view, in order to compute the considered probability, there is no need to know the complete enumeration of the sequence  $(|Av_m(\sigma)|)_m$. Indeed, the sum in (\ref{buffetto}) converges because so does $\sum_{m=1}^{\infty} |Av_m(\sigma)|/(m-1)!$. For any $2 \leq k \leq n$ we thus define
\begin{equation}\label{tilda}
w_{\sigma,n,k} = \sum_{m=1}^{n-k+1} \frac{{{n-m-1}\choose{k-2}}}{{{n-2}\choose{k-2}}} \times \frac{|Av_m(\sigma)|}{(m-1)!}
\end{equation} 
and we have the following
\begin{prop}\label{pviv}
For a random permutation $\pi$ of size $n$ and $2 \leq k \leq n$ we have 
\begin{equation}\label{vivaldi}
\text{Prob}(\pi \notin Av_n(\sigma;k)) \simeq w_{\sigma,n,k} \cdot \frac{k(k-1)}{n(n-1)},
\end{equation}
where $w_{\sigma,n,k}$ converges to a constant (\ref{tilda}) when $n$ is large enough.
\end{prop}

\bigskip

\begin{figure}
\begin{center}
\includegraphics*[scale=.85,trim=0 0 0 0]{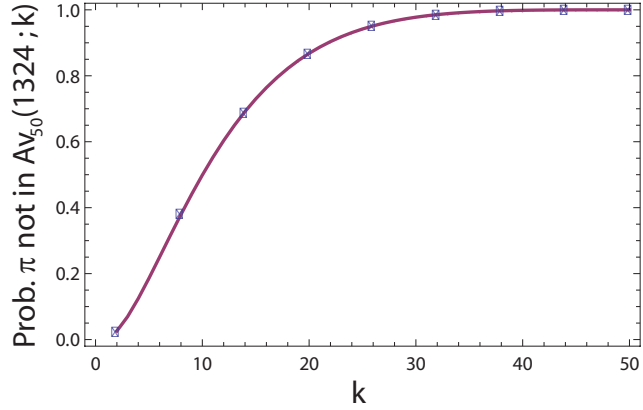}
\end{center}
\caption{Plot of $\text{Prob}(\pi \notin Av_{50}(1324;k))$ for $2\leq k \leq 50$. Points are obtained experimentally via random generation. The line is drawn according to (\ref{vivaldi}) using the first $20$ terms of the sequence $|Av_m(1324)|$.}\label{ro}
\end{figure}

For instance, take $\sigma=1324$ and $n=50$. Consider the first $20$ terms of the sequence $|Av_n(1324)|$ (\cite{sloane} seq. A061552). We approximate $w_{1324,50,k}$ as $$w_{1324,50,k} \simeq \sum_{m=1}^{\min(20,50-k+1)} \frac{{{50-m-1}\choose{k-2}}}{{{50-2}\choose{k-2}}} \times \frac{|Av_m(\sigma)|}{(m-1)!}.$$ Then we can test the estimate given in (\ref{vivaldi}) with data from simulations for $n=50$ (middle line in Fig.~\ref{vsk}). The result is depicted in Fig.~\ref{ro}.

\section*{Acknowledgement}
This work was financially supported by grant DFG-SPP1590 from the German Research Foundation to TW.

\end{document}